\documentclass[a4paper,11pt,centertags,psamsfonts]{amsart}
\usepackage{amssymb}
\usepackage{times}
\usepackage[mathcal]{euscript}

\copyrightinfo{2006}{V.~Kostrykin, K.~A.~Makarov}

\renewcommand{\epsilon}{\varepsilon}

\newcommand{\R}{\mathbb{R}}

\newcommand{\C}{\mathbb{C}}

\newcommand{\N}{\mathbb{N}}

\newcommand{\EE}{\mathsf{E}}

\newcommand{\cL}{{\mathcal L}}

\newcommand{\cP}{{\mathcal P}}

\renewcommand{\Re}{{\ensuremath{\mathrm{Re}}}}

\newcommand{\tr}{\mathrm{tr}}

\DeclareMathOperator{\Ran}{\mathrm{Ran}}
\DeclareMathOperator{\Ker}{\mathrm{Ker}}

\newcommand{\ii}{\mathrm{i}}
\newcommand{\e}{\mathrm{e}}

\newtheorem{theorem}{Theorem}{\bf}{\it}
{\bf}{\it}
{\bf}{\it}
{\bf}{\it}
{\bf}{\it}
\theoremstyle{definition}
{\bf}{\it}
{\it}{\rm}

\newtheorem{remark}[theorem]{Remark}{\it}{\rm}

\title[On Krein's Example]{On Krein's Example}

\author[V. Kostrykin]{Vadim Kostrykin}
\address{V.~Kostrykin, Institut f\"{u}r Mathematik, Technische Universit\"{a}t Clausthal,
Erzstra{\ss}e 1, D-38678 Clausthal-Zellerfeld, Germany}
\email{kostrykin@math.tu-clausthal.de, kostrykin@t-online.de}

\author[K. A. Makarov]{Konstantin A.~Makarov}
\address{K.~A.~Makarov, Department of Mathematics, University of
Missouri, Co\-lum\-bia, MO 65211, USA}
\email{makarov@math.missouri.edu}

\keywords{Spectral shift function, Hankel operators, absolutely continuous
spectrum.}

\subjclass[2000]{Primary 47B35; Secondary 47A55, 45P05}

\begin{document}

\begin{abstract}
In his 1953 paper [Matem.~Sbornik \textbf{33} (1953), 597 -- 626] Mark
Krein presented an example of a symmetric rank one perturbation of a
self-adjoint operator such that for all values of the spectral parameter in
the interior of the spectrum, the difference of the corresponding spectral
projections is not trace class. In the present note it is shown that in the
case in question  this difference has simple Lebesgue spectrum filling in
the interval $[-1,1]$ and, therefore, the pair of the spectral projections
is generic in the sense of Halmos but not Fredholm.
\end{abstract}

\maketitle

The spectral shift function plays a very important role in perturbation
theory for self-adjoint operators. It was introduced in a special case by
I.~Lifshitz \cite{Lifshitz:52} and in the general case (in the framework of
trace class perturbations) by M.~Krein in his celebrated 1953 paper
\cite{Krein}. He showed that for a pair of self-adjoint not necessarily
bounded operators $A_0$ and $A_1$ such that their difference $A_1-A_0$ is
trace class there exists a unique function $\xi\in L^1(\R)$ satisfying the
trace formula
\begin{equation}\label{trace:formula}
\tr(\varphi(A_1)-\varphi(A_0)) = \int_\R \varphi^\prime(x) \xi(x) dx
\end{equation}
whenever $\varphi$ belongs to a class of admissible functions.
Surprisingly, the condition $\varphi\in C^1$ alone does not imply the validity of the
trace formula \eqref{trace:formula} (see \cite{Peller}).

If one goes ahead and  formally puts $\varphi=\chi_\mu$ in \eqref{trace:formula} with
$\chi_\mu$ a characteristic function of the interval $\delta_\mu :=
(-\infty,\mu)$, one would arrive at the ``naive Lifshitz formula''
\begin{equation}\label{naive}
\xi(\mu) = \tr(\EE_{A_1}(\delta_\mu)-\EE_{A_0}(\delta_\mu)),
\end{equation}
where $\EE_{A_j}(\delta_\mu)$ denotes the spectral projection for the
operator $A_j$, $j=0, 1$, associated with the semi-infinite interval $\delta_\mu$.
Note, that representation \eqref{naive} holds for all $\mu\in \Delta$ whenever
$\Delta$ is a joint spectral gap for the operators $A_0$ and $A_1$.

In general the right hand side of \eqref{naive}
does not make sense: In Section 6 of \cite{Krein}
Krein presented an example of two self-adjoint operators such that their
difference $A_1-A_0$ is of rank one but
$\EE_{A_1}(\delta_\mu)-\EE_{A_0}(\delta_\mu)$ is not trace class and, hence, the
trace in the r.h.s.\ of \eqref{naive} is ill-defined (we refer the reader  to
the review \cite{Birman:Pushnitski} for the further discussion).

In this example Krein introduces a bounded integral
operator $A_0$ in $L^2(0, \infty)$
with kernel given by
\begin{equation*}
A_0(x,y)=\begin{cases}
                \sinh x\, \e^{-y}, & x\leq y\\
                \sinh y\, \e^{-x}, & x\geq y
                \end{cases}
\end{equation*}
 and its rank one perturbation $A_1$ with kernel
\begin{equation*}
A_1(x, y)=A_0(x, y)+ \e^{-x} \e^{-y}
\end{equation*}
so that
 \begin{equation*}
A_1(x,y)=\begin{cases}
                \cosh x\, \e^{-y}, & x\leq y\\
                \cosh y\, \e^{-x}, & x\geq y
                \end{cases}.
\end{equation*}
In fact $A_j$, $j=0,1$, are resolvents of the Dirichlet ($j=0$)
and Neumann ($j=1$) one-dimensional Laplacian $-d^2/dx^2 $ at the
spectral point $-1$. Obviously, the operators $A_j$, $j=0,1$, have
simple, purely absolutely continuous spectrum filling in the
interval $[0,1]$.

Krein shows that the difference
\begin{equation}\label{K:mu}
K_\mu := \EE_{A_1}(\delta_\mu)-\EE_{A_0}(\delta_\mu),
\quad 0<\mu < 1,
\end{equation}
of the spectral projections for the operators $A_1 $ and $A_0$ is
the integral operator with the kernel $k_\mu(x+y)$, where
\begin{equation*}
k_\mu (x)=\frac{2}{\pi}\frac{\sin\sqrt{\lambda(\mu)}\,x}{x}
\end{equation*}
with
\begin{equation*}
\lambda(\mu)=\frac{1}{\mu}-1,\quad 0<\mu <1.
\end{equation*}

As noted by Krein in \cite{Krein} the  operator  $K_\mu$, $0<\mu
<1$, is not Hilbert-Schmidt, since for any $\lambda>0$
\begin{equation*}
\int_0^\infty \int_0^\infty \frac{\sin^2\sqrt{\lambda} (x+y)}{(x+y)^2} dx
dy= \int_0^\infty \frac{\sin^2 u}{u}du=\infty.
\end{equation*}

Employing the theory of Hankel operators \cite{Howland:85} it not
hard to show that the operator $K_\mu$, $0<\mu<1$, has an
absolutely continuous spectrum $[-1,1]$ of uniform multiplicity
one. However, to the best of our knowledge there are no general
criteria guaranteeing the absence of singular spectrum for this
operator. In particular, the related results of
\cite{Howland:92a}, \cite{Howland:92b} are not applicable to the
case in question.

In the present note we perform a detailed spectral analysis of the
operator $K_\mu$. Our result is given by the following

\begin{theorem}\label{Thm:1}
For any $0<\mu<1$ the operator $K_\mu$ defined in \eqref{K:mu}
has a simple purely absolutely
continuous spectrum filling in the interval $[-1,1]$.
\end{theorem}

\begin{proof}
Observe that $K_\mu$, $0<\mu<1$, is unitary equivalent
to $K_{\frac12}$ with the unitary equivalence given by the
scaling transformation
\begin{equation*}
U_\lambda: f(x)\mapsto \lambda^{1/4}
f(\sqrt{\lambda} x), \quad f\in L^2(0, \infty).
\end{equation*}
That is,
\begin{equation*}
K_\mu=U_{\lambda(\mu)} K_{\frac12}
U_{\lambda(\mu)}^\ast.
\end{equation*}

Therefore, it suffices to prove the theorem for the operator
$K_{\frac12}$. For brevity we set $K:=K_\frac12$. Since
\begin{equation*}
\frac{2}{\pi}\frac{\sin\, x}{x}=\frac{1}{2\pi} \int_{-1}^1 2\, \e^{-\ii tx}
dt,
\end{equation*}
the operator $K$ is a Hankel integral operator
with discontinuous symbol
\begin{equation*}
s(t) = \begin{cases} 2, & t\in[-1,1] \\ 0, & t\notin[-1,1] \end{cases}.
\end{equation*}

In the Hilbert space  $\ell_+^2$ of complex square-summable one-sided
sequences $x=\{x_0, x_1, \ldots \}$ introduce
the Hankel operator $H(\phi)$
with symbol $\phi$
given by
\begin{equation*}
\phi(z) := s\left(\ii\frac{1-z}{1+z}\right),\quad |z|=1,\quad z\neq -1,
\end{equation*}
which is the characteristic function of the set $\{z\in \C | \,|z|=1,\,
\Re\, z\geq 0\}$. That is,
\begin{equation*}
(H(\phi)x)_n = \sum_{k=0}^\infty c_{n+k+1} x_k,\quad x=\{x_0, x_1, ... \}\in
\ell_+^2.
\end{equation*}
where $c_k$ are Fourier coefficients of the function $\phi$,
\begin{equation*}
c_k = \frac{1}{2\pi} \int_0^{2\pi} \e^{\ii k \theta} \phi(\e^{\ii\theta})
d\theta = \frac{2}{\pi k} \sin(\pi k/2),\qquad k\in\N.
\end{equation*}
The operators $K$ and $H(\phi)$ are unitarily equivalent (see, e.g.,
\cite[p.~14]{Power}).
For any $p\notin\N$ let $H_{p}$ and $\widetilde{H}_{p}$ be the $p$-shifted
Hilbert matrices on $\ell^2_+$,
\begin{equation*}
(H_{p} x)_n = \sum_{k=0}^\infty (n+k+1-p)^{-1} x_k,\qquad
(\widetilde{H}_{p} x)_n = \sum_{k=0}^\infty (-1)^{n+k}(n+k+1-p)^{-1} x_k.
\end{equation*}
In particular, $H_0$ is the standard Hilbert matrix.

Denote by $\cL_+$ (respectively $\cL_-$) the set of all sequences in
$\ell_+^2$ with vanishing odd (respectively even) elements, that is,
\begin{equation*}
\begin{split}
\cL_+ & = \{x\in\ell_+^2| x_{2k+1}=0\quad \forall k\in\N_0\},\\
\cL_- & = \{x\in\ell_+^2| x_{2k}=0\quad \forall k\in\N_0\}.
\end{split}
\end{equation*}
Let $\cP_\pm$ be the orthogonal projections in $\ell^2_+$ onto $\cL_\pm$,
respectively.

Observe that $c_k=0$ for all $k\in 2\N$, $c_k=\frac{2}{\pi k}$ for all
$k\in 4\N_0+1$, and $c_k=-\frac{2}{\pi k}$ for all $k\in 4\N_0+3$. This
immediately implies that $\cP_+ H(\phi)\cP_- = 0$ and $\cP_- H(\phi)\cP_+ =
0$. Furthermore, for any $x\in\cL_+$ we have
\begin{equation*}
(\cP_+ H(\phi)\cP_+ x)_{2n} = \frac{2}{\pi}\sum_{k=0}^\infty
\frac{(-1)^{n+k}}{2(n+k+1/2)} x_{2k}.
\end{equation*}
Thus, $\cP_+ H(\phi)\cP_+$ is unitarily equivalent to the operator
$\pi^{-1}\widetilde{H}_{1/2}$ on $\ell_+^2$. Similarly, for any $x\in\cL_-$
we have
\begin{equation*}
(\cP_- H(\phi)\cP_- x)_{2n+1} = -\frac{2}{\pi}\sum_{k=0}^\infty
\frac{(-1)^{n+k}}{2(n+k+3/2)} x_{2k+1},
\end{equation*}
which implies that $\cP_- H(\phi)\cP_-$ is unitarily equivalent to
$(-\pi^{-1}\widetilde{H}_{-1/2})$ on $\ell_+^2$.

Since for any $p\notin\N$ the operators $H_{p}$ and $\widetilde{H}_{p}$ are
unitarily equivalent, we have that the Hankel operator $H(\phi)$ is
unitarily equivalent to the orthogonal sum $\pi^{-1}H_{1/2}\oplus
(-\pi^{-1}H_{-1/2})$ with
respect to the orthogonal decomposition $\ell_+^2\oplus \ell_+^2$.

By a result of Rosenblum \cite[Theorem 5]{Rosenblum:2} for any $p\leq 1/2$
the spectrum of the operator operator $H_p$ is purely absolutely
continuous, has uniform multiplicity $1$, and fills the interval $[0,\pi]$.
This completes the proof of the theorem.
\end{proof}

We conclude our note by a remark of geometric character.

\begin{remark}
Since  for any $\mu\in (0, 1)$ the essential spectrum of  $K_\mu$
fills in the interval $[-1,1]$, one concludes that the pair
$(\EE_{A_1}(\delta_\mu),\EE_{A_0}(\delta_\mu))$ of the spectral
projections is not Fredholm in the sense of \cite{Avron}. Moreover,
the absence of the discrete spectrum of $K_\mu$ implies that the
subspaces $\Ran\EE_{A_1}(\delta_\mu)$ and
$\Ran\EE_{A_0}(\delta_\mu)$ are in generic position in the sense
of Halmos \cite{Halmos}. In particular,
\begin{equation*}
 \Ker( \EE_{A_1}(\delta_\mu)-\EE_{A_0}(\delta_\mu)\pm I)=\{0\}
\end{equation*}
for all $\mu\in(0,1)$.
\end{remark}


\end{document}